\documentclass[12pt]{amsart}
\newtheorem{theorem}{Theorem}
\newtheorem{corollary}{Corollary}

\theoremstyle{definition}

\newtheorem*{example}{Example}

\theoremstyle{remark}
\newtheorem*{remark}{Remark}

\newcommand{\Q}{\mathbf{Q}}

\newcommand{\C}{\mathbf{C}}

\newcommand{\A}{\mathbf{A}}

\begin{document}

\title[On strong multiplicity one]{On strong multiplicity one for automorphic representations}

\author{C. S. Rajan}

\address{School of Mathematics, Tata Institute of Fundamental 
Research, Homi Bhabha Road, Bombay - 400 005, INDIA.}

\email{rajan@math.tifr.res.in}


\begin{abstract}
We extend the strong multiplicity one theorem of Jacquet,
Piatetski-Shapiro and Shalika. Let $\pi$ be a unitary, cuspidal,
automorphic representation of $GL_n(\A_K)$. Let $S$ be a set of finite
places of $K$, such that the sum $\sum_{v\in S}Nv^{-2/(n^2+1)}$ is
convergent.  Then $\pi$ is uniquely determined by the collection of
the local components $\{\pi_v\mid v\not\in S, ~v ~\text{finite}\}$ of
$\pi$.  Combining this theorem with base change, it is possible to
consider sets $S$ of positive density, having appropriate splitting
behavior with respect to a solvable extension $L$  of $K$, and where $\pi$ is
determined  upto twisting by a character of the Galois group of $L$
over $K$.  
\end{abstract}

\maketitle

\section{Introduction} 
Let $K$ be a number field with ring of adeles $\A_K$.  Let ${\mathcal
A}(n,K)$ denote the set of isomorphism classes of irreducible,
unitary, cuspidal  automorphic representations of $GL_n(\A_K)$. Given
$\pi\in {\mathcal A}(n,K)$, $\pi$ can be expressed as a restricted tensor
product $`\otimes'_{v\in \Sigma_K}\pi_v$, where $v$ runs over all the
$\Sigma_K$ places of $K$,
and $\pi_v$ is an irreducible, unitary representation of
$GL_n(K_v)$. At  almost all finite places $v$, $\pi_v$ is
unramified, and let  $A_v(\pi)$ denote the conjugacy class in
$GL(n, \C)$ determined by the Langlands-Satake parameter associated to
$\pi_v$. 

Let $S$ be a finite set of places of $K$.  The strong multiplicity one
theorem of Jacquet, Piatetski-Shapiro and Shalika \cite{JS}, asserts
that any $\pi\in {\mathcal A}(n,K)$ is uniquely determined by the
collection of the local components $\{\pi_v\mid v\in \Sigma_K- S\}$ of
$\pi$. In view of the possible applications to the functoriality
principle of Langlands, it is desirable to have refined versions of
the strong multiplicity one theorem. In this note, we give an
analytical proof of a stronger version of the strong mulitplicity one
theorem, based on the analytical properties of the Rankin-Selberg
$L$-functions.

Let ${\mathcal O}_K$ denote the ring of integers of $K$. 
 If $v$ is a finite place of $K$, then we
denote the corresponding prime ideal by $P_v$, and by $Nv$ the norm
of $v$ to be the number of elements in the residue field ${\mathcal
 O}_K/P_v$.  
Let $S$ be a set of  finite places of $K$  
satisfying the following convergence condition: 
\begin{equation}\label{convcriterion}
 \sum_{v\in S}Nv^{-\frac{2}{n^2+1}}< \infty.  
\end{equation}
Our main theorem is as follows: 
\begin{theorem}  \label{maintheorem}
Let  $\pi_1$ and $\pi_2$ be irreducible, unitary, cuspidal automorphic
representations of $GL_n(\A_K)$.
 Let $S$ be a set of finite places of $K$   satisfying the
convergence criterion given in \eqref{convcriterion}.  Suppose
that at almost all but  finitely many  places $v\not\in S$, the local
components of $\pi_1$ and $\pi_2$ are isomorphic:  
\[ \pi_{1, v}\simeq \pi_{2,v}.\]
Then $\pi_1\simeq \pi_2.$
\end{theorem}

An advantage with this formulation,  is that the ignorant
set of places $S$ can be infinite, and in fact the convergence
criterion allows us to handle sets $S$ of positive density. 

\begin{example}
Let $K/\Q$  be a cyclic
extension of prime degree greater than $(n^2+1)/2$ over $\Q$, and take
$S$ to be the collection of unramified finite places in $K$, which are
not of degree one over $\Q$. 

More generally, let $K/\Q$ be a Galois extension with Galois group
$G$.  Let $\sigma_v$ denote the Frobenius element at an unramified
finite place $v$ of $K$ over $\Q$. Then the collection $S$ of finite
unramified places $v$ of $K$, such that the order of $\sigma_v$ is
greater than $(n^2+1)/2$ satisfies the convergence condition
(\ref{convcriterion}).
\end{example}

The above theorem,  combined with the results on solvable base change
of automorphic representations established by Langlands, Arthur and
Clozel \cite{AC}, and the characterisation of the fibres of solvable
base change \cite{R3}, allow us to handle sets $S$ of
positive density, consisting of places of $K$ 
 having suitable splitting properties with respect to
solvable extensions of $K$:
\begin{corollary} Suppose $\pi_1$ and $\pi_2$ are irreducible,
unitary, cuspidal automorphic representations of $GL_n(\A_K)$. 
  Let $L/K$ be a solvable
 extension, and $S$ be a set of finite places in $K$, unramified in
 $L$, such that for any place $w$ of $L$ dividing a place $v$ in
 $S$ of residue characteristic $p$, the norm $Nw$ is greater than
 $p^{(n^2+1)/2}$.  Suppose that at almost all but finitely many  places
 $v$ of $K$ not belonging to $S$, the local components  of
 $\pi_1$ and $\pi_2$ are isomorphic:
\[ \pi_{1, v}=\pi_{2,v}.\]
Assume further that the base change of $\pi_1$ and $\pi_2$ to
$GL_n(\A_L)$ is cuspidal. Then there is an idele class character $\chi$
of $K$, corresponding via class field theory to a character of the
Galois group of the extension $L/K$, such that
\[\pi_1\simeq \pi_2\otimes \chi.\]
\end{corollary}
\begin{proof} Let $BC_{L/K}$ be the base change map defined by Arthur
and Clozel (\cite{AC}), from ${\mathcal A}(n,K)$ to ${\mathcal
A}(n,L)$. By Theorem \ref{maintheorem}, it follows from the hypothesis
of the corollary, that $BC_{L/K}(\pi_1)\simeq BC_{L/K}(\pi_2)$. Since
the base change representations are cuspidal, the corollary  follows from the
description of the fibres of solvable base change given by Theorem 2
of \cite{R3}. 
\end{proof}

A conjecture of Ramakrishnan \cite{DR2}, asserts that if the Dirichlet
density of the set of places $S$ of $K$ is less than $1/2n^2$, then
the set of local components $\pi_v$ for $v$ not in $S$, determines
$\pi$ uniquely. Ramakrishnan in \cite{DR1}, proved this conjecture in
the case $n=2$, using the automorphicity of the symmetric square
lifting for $GL(2)$. In \cite{R2}, Ramakrishnan's conjecture is
verified for a particular class of automorphic representations of
$GL_n(\A_K)$, and also that Ramanujan's conjecture (and not a weak
Ramanujan conjecture, as wrongly stated there) implies Ramakrishnan's
conjecture.

In the $l$-adic context, apart from Ramakrishnan's conjecture, much
more is known about the nature of strong multiplicity one results
\cite{R1}. For example, suppose $\rho_1, ~\rho_2:G_K\to GL_n(\Q_l)$
are absolutely irreducible $l$-adic representations of the absolute
Galois group $G_K$ of a global field $K$, unramified outside a finite
set of places. Suppose that the algebraic envelope of the image of one
of the representations is connected, and that the traces of the
Frobenius elements agree at a set of unramified finite places of
positive density of $K$. Then there is a character $\chi: G_K\to
\Q_l^*$  of finite order, such that $\rho_2\simeq \rho_1\otimes \chi$. 

Analogously, it can  be expected \cite{R2}, that if $\pi_1$ and
$\pi_2$ are cuspidal, automorphic representations of `general type' on
$GL_n(\A_K)$ (when $n=2$, general type is the same as non
automorphically induced from a Hecke character of a quadratic
extension of $K$), 
 such that $a_v(\pi_1)=a_v(\pi_2)$ at a set of unramified
finite places of positive density, then there exists an idele-class
character $\chi$ of $K$, such that $\pi_1\simeq \pi_2\otimes
\chi$. Here for an unramified finite place $v$ of an automorphic
representation $\pi$, the Dirichlet coefficient  $a_v(\pi)$ of $\pi$
at $v$,  is  the trace of $A_v(\pi)$.  The above results represent a
small step towards this general direction.

\section{Proof}
The method of proof is analogous to the proof of Ramakrishnan's  
conjecture for Artin (finite)  representations  of  the absolute
Galois group:
let $G$ be a finite group, and let $\rho_1,~\rho_2$ be
inequivalent representations of $G$ to $GL(n,\C)$. Then 
\[ \#\{ g\in G\mid Tr(\rho_1(g))=Tr(\rho_1(g))\}
\leq (1-1/2n^2)|G|.
\]
The proof of this result follows from considering the equality, 
\[
\frac{1}{|G|}\sum_{g\in G}|Tr(\rho_1(g))-Tr(\rho_2(g))|^2\geq 2, 
\]
and using the fact that the character value at an element $g$ in $G$,
is a sum of $n$ roots of unity.

In analogy, we consider the $L$-function associated to an automorphic
representation,  as the `character' of the representation. Then the
convolution product of characters becomes the Rankin-Selberg
$L$-function $L(s,\pi_1\times\pi_2)$ associated to a pair of unitary,
cuspidal automorphic representations $\pi_1$ and $\pi_2$  of
$GL_n(\A_K)$.  The analytic properties of the convolution $L$-function
have been studied extensively by Jacquet, Piatetskii-Shapiro, Shalika,
Shahidi and Waldspurger in a series of papers \cite{JS}, \cite{JPSh},
\cite{Sh}, \cite{MoW}.

The local factors of the Rankin-Selberg $L$-function 
$L(s, \pi_1\times \pi_2)$ at an unramified finite place $v$ of $\pi_1$ and
$\pi_2$  is defined as, 
\begin{equation}\label{unramifiedfactor}
L(s, \pi_{1,v} \times \pi_{2,v})={\rm 
det}(1-A_v(\pi_{1})\otimes A_v(\pi_{2})Nv^{-s})^{-1} .
\end{equation}
Let $T$ be a set of  places of $K$ containing the archimedean places
of $K$.  Define the Dirichlet series
\[ L_T(s, \pi_1 \times \pi_2)=\prod_{v\not\in T}L(s,
\pi_{1,v} \times \pi_{2,v}). \]
If $T$ is finite, then $L_T(s, \pi_1 \times \pi_2)$ satisfies the following
properties (\cite{JS}, \cite{JPSh}, \cite{Sh},
\cite{MoW}):
\begin{itemize}
\item The Dirichlet series $ L_T(s, \pi_1 \times \pi_2)$ is absolutely convergent in the half
plane $\text{Re}(s)>1$. 
\item The function $L_T(s, \pi_1\times \pi_2)$ admits a meromorphic
continuation to $1$. Further, if $L_T(s, \pi_1\times \pi_2)$ is
holomorphic at $s=1$, then it is  non-vanishing at $s=1$.   
\item Let
$\tilde{\pi}$ denote the contragredient representation of an
automorphic representation $\pi$.  The function $L_T(s, \pi_1\times
\pi_2)$ has a simple pole at $s=1$ if and only if $\pi_1\simeq \tilde{\pi}_2$.
\end{itemize}
Consider the following Dirichlet series
\[ L_T(s) =\frac{L_T(s, \pi_1\times\tilde{\pi}_1)L_T(s,
 \pi_2\times\tilde{\pi}_2)} {L_T(s, \pi_1\times\tilde{\pi}_2)L_T(s,
\pi_2\times\tilde{\pi}_1)}\] 
It follows from the properties of
Rankin-Selberg convolutions listed above,  that if $T$ is finite, then
$L_T(s)$ has a pole of order $2$, provided  $\pi_1\not\simeq  \pi_2$. 

Now let $T$ be a finite set of places of $K$ containing the
archimedean places of $K$, the ramified places of $\pi_1$ and $\pi_2$,
and the finite number of places $w$ not in $S$, as in the hypothesis
of the theorem, where it is not known that the local $\pi_{1,w}$ and
$\pi_{2,w}$ are isomorphic.  Let $S', T'$ denote respectively the
complement of $S$ and $T$ in the set of all places of $K$. We have the
hypothesis of the theorem that the local components $\pi_{1,v}$ and
$\pi_{2,v}$ are isomorphic at all places of $K$ outside $T$ and $S$.
Hence we obtain,
\begin{equation}\label{theproduct}
L_T(s)=L_{T\cup S'}(s)=\prod_{v\in S\cap T'}
\frac{L(s, \pi_{1,v}\times\tilde{\pi}_{1,v})L(s,
 \pi_{2,v}\times\tilde{\pi}_{2,v})} {L(s, \pi_{1,v}\times\tilde{\pi}_{2,v})L(s,
\pi_{2,v}\times\tilde{\pi}_{1,v})}
\end{equation}

At an unramified finite place $v$ of an automorphic representation
$\pi$, with associated Langlands-Satake parameter $A_v(\pi)$, define
for any natural number $k$, 
\[ a_{v,k}(\pi)=\sum_{i=1}^n\alpha_{v,i}(\pi)^k,\]
where $\alpha_{v,1},\cdots, \alpha_{v,n}$ are the eigenvalues of
matrix representative of $A_v(\pi)$.  Expanding ${\rm log}~L_T(s)$ in
terms of the Euler product expansion given by \eqref{theproduct}, and
by the definition of the local $L$-factors as given by
\eqref{unramifiedfactor}, we obtain in the region  $\sigma =\text{Re}(s)>1$
\begin{equation}\label{logfunction}
 {\rm log}~L_T(s)= \sum_{v\in S\cap T'}\sum_{k=1}^{\infty}
 |a_{v,k}(\pi_1)-a_{v,k}(\pi_2)|^2Nv^{-ks}k^{-1}.
\end{equation}
We recall the estimate for an eigenvalue $\alpha_v$ of $A_v(\pi)$
proved by Luo, Rudnick and Sarnak \cite{LRS} for an unramified finite
place $v$ of a unitary, automorphic representation $\pi$ of
$GL_n(\A_K)$:
\begin{equation}\label{LRS}
 |\alpha_v|\leq Nv^{\frac{1}{2}-\frac{1}{n^2+1}}.
\end{equation} 
It follows from this estimate, that for $\sigma>1$  we have a majorization 
\[\left|{\rm log}~L_T(s)\right|\leq  4n^2 \sum_{v\in S\cap T'}
\sum_{k=1}^{\infty}k^{-1}Nv^{-k\sigma+k-2k/(n^2+1)}. \]
The  convergence
condition satisfied by $S$ as given in  \eqref{convcriterion},
ensures that the right hand expression has a finite limit as
$s\rightarrow 1^+$. But this implies that $\pi_1\simeq \pi_2$, and
that proves the theorem. 

\begin{remark} The theorem can be extended to pairs of isobaric,
unitary automorphic representations $\pi_1$ and $\pi_2$ of $GL_n(\A_K)$. 
\end{remark}

\begin{remark} It is possible to use the method to compare two
isobaric, unitary automorphic representations $\pi_1$ and $\pi_2$ on
$Gl_n(\A_K)$ and $GL_m(\A_K)$ respectively ($n\geq m$). Instead of
requiring the equality of local components, we impose the hypothesis
\[ a_{v,k}(\pi_1)=a_{v,k}(\pi_2) \quad k=1,\cdots, [(n^2+1)/2].\]
The Luo-Rudnick-Sarnak estimates \eqref{LRS}, ensure that the
contribution of the sum over the places $v\in S'$ to ${\rm
log}~L_T(s)$ is $O(1)$ as $s\rightarrow 1^+$, and the rest of the proof
goes through. 

If we assume further the Ramanujan conjectures, then the method of
proof allows us to show that if the Dirichlet density of $S$ is less than
$1/2n^2$, then the collection of Dirichlet coefficients (rather than
the local components) $\{a_v(\pi)\mid v\not \in S\}$ uniquely
determines $\pi$. 
\end{remark}

\begin{remark} We give a brief indication of the proof of
a slight strenthening of Ramakrishnan's theorem \cite{DR1}, that if
the Dirichlet density 
 $d:= \lim_{s\rightarrow 1^+} -{\sum_{v\in S}Nv^{-s}}/{\rm
log}~(s-1)$ of $S$ is at most $1/8$, then the collection of Dirichlet
coefficients $\{a_v(\pi)\mid v\not\in S\}$ of $\pi$ (rather than local
components) uniquely determines $\pi$, where $\pi$ is a cuspidal
automorphic of $GL_2(\A_K)$.  

Using the estimate 
$|\alpha_v(\pi)|\leq Nv^{1/5}$ proved in \cite{LRS} (even the weaker
estimate $Nv^{1/4-\epsilon}$ for some positive $\epsilon$ will
suffice), we see that the
hypothesis implies, 
\begin{equation}\label{dinakar1}
\begin{split}  {\rm log}~ &  L_T(s)  =\sum_{v\in S\cap T'}
 |a_{v}(\pi_1)-a_{v}(\pi_2)|^2Nv^{-s} +O(1)\quad
\text{as} ~s\rightarrow 1^+\\
& \leq 2\sum_{v\in S\cap T'}
 (|a_{v}(\pi_1)|^2+
 |a_{v}(\pi_2)|^2)Nv^{-s}+O(1)\quad
\text{as} ~s\rightarrow 1^+.
\end{split}
\end{equation}  
For $\pi\in {\mathcal A}(2,K)$, let $\text{Ad}(\pi)$ denote the
adjoint lift of $\pi$, an automorphic representation of $GL(3,\A_K)$
constructed by Gelbart and Jacquet \cite{GJ}.
  If $\text{Ad}(\pi)$ is not cuspidal,
then it is automorphically induced and we have the estimate $
|a_{v}(\pi)|^2\leq 4$. Hence we obtain, 
\begin{equation}\label{dihedral}
-\lim_{s\rightarrow 1}\frac{\sum_{v\in S\cap T'}
 |a_{v}(\pi)|^2Nv^{-s}}{{\rm log}~(s-1)}=4d.
\end{equation}
If  $\text{Ad}(\pi)$ is cuspidal, then consider
\[L^4(s):= L_T(s, \pi\times \tilde{\pi}\times \pi\times \tilde{\pi})
=L_T(s, \text{Ad}(\pi)\times \text{Ad}(\pi))L_T(s,
\text{Ad}(\pi))^2\zeta_{K,T}(s).  \] 
From the properties of the
Rankin-Selberg $L$-functions we obtain
\[{\rm log}~L^4(s)= -2{\rm
log}(s-1)+f(s)=\sum_{v\not\in T}\frac{|a_v(\pi)|^4}{Nv^s} +O(1) \quad
\text{as} ~s\rightarrow 1^+,\] 
 where $f(s)$ is holomorphic at $s=1$.  By  
Cauchy-Schwarz, we get
\begin{equation}\label{ad-cuspidal}
-\lim_{s\rightarrow 1}\frac{\sum_{v\in S\cap T'}
 |a_{v}(\pi)|^2Nv^{-s}}{{\rm log}~(s-1)}=\sqrt{2d}.
\end{equation}
From \eqref{dinakar1}, \eqref{dihedral} and \eqref{ad-cuspidal}, we
obtain, 
\[ 2= \lim_{s\rightarrow 1} -\frac{{\rm log}~L_T(s)}{{\rm log}~(s-1)} \leq
\begin{cases} 4\sqrt{2d}& \text{if $\text{Ad}(\pi_1)$ and
$\text{Ad}(\pi_2)$ are cuspidal},\\
 2\sqrt{2d}+8d &  \text{if only  one of $\text{Ad}(\pi_1)$ or
$\text{Ad}(\pi_2)$ is cuspidal}, \\
 16d & \text{if $\text{Ad}(\pi_1)$ and
$\text{Ad}(\pi_2)$ are not cuspidal}.
\end{cases}
  \]
In all the three cases, this leads to  a contradiction if $d<1/8$. 
This proves Ramakrishnan's theorem. 

It is not clear whether the methods contained in this paper, will be
sufficient to prove that $S$ can be taken to be a set of places of
density at most $1/2$,  provided that  both $\text{Ad}(\pi_1)$ and
$\text{Ad}(\pi_2)$ are cuspidal. Such a result can be expected in
analogy with the results in the $l$-adic context (and hence valid also
for holomorphic modular forms) obtained in \cite{R1}. 

Our original approach to the results of this paper was to use the
Cauchy-Schwarz inequality,  and to obtain a region of convergence to the
left of $2$ for the sum $\sum_{v}|a_v(\pi)|^4Nv^{-s}$. Using the
methods contained in \cite{DI}, it is possible to obtain a region of
convergence to the left of $2$, provided $K=\Q$. But over number
fields, it is not clear how to obtain such a result. However, even
over $\Q$, these results do not yield stronger results than our main
theorem towards strong multiplicity one.   
\end{remark}

\end{document}